# STURM-LIOUVILLE PROBLEMİNİN REZOLVENT OPERATÖRÜ VE ÖZFONKSİYONLARI

## Erdoğan ŞEN[1], Oktay MUKHTAROV[2], Kamil ORUÇOĞLU[3]

[1] *Namık Kemal Üniversitesi, Fen Edebiyat Fakültesi, Matematik Bölümü, 59030, Tekirdağ*
[2] *Gaziosmanpaşa Üniversitesi, Fen Fakültesi, Matematik Bölümü, 60250, Tokat*
[1,3] *İTÜ Fen-Edebiyat Fakültesi Matematik Mühendisliği Bölümü, 34469, Maslak, İstanbul*

**Abstract**

*In this work we investigate the resolvent operator and completeness of eigenfunctions of a Sturm-Liouville problem with discontinuities at two points. The problem contains an eigenparameter in the one of boundary conditions. For operator-theoretic formulation of the considered problem we define an equivalent inner product in the Hilbert space $L_2[-1,1]\oplus$ and suitable self-adjoint lineer operator in it.*

**Özet**

*Bu çalışmada iki noktada süreksiz olan ve sınır koşullarının birinde özdeğer parametresi içeren Sturm-Liouville probleminin rezolvent operatörünü ve özfonksiyonlarının tamlığını inceledik. İncelenen problemin operatör-kuramsal yazılımı için $L_2[-1,1]\oplus$ Hilbert uzayında yeni eşdeğer iç çarpım ve uygun kendine eşlenik lineer operatör tanımlandı.*

***Anahtar Kelimeler:*** *Sturm-Liouville problemi, özdeğer, özfonksiyon, rezolvent operatör, geçiş koşulları.*

## GİRİŞ

Muhtarov ve diğerleri [7]'de sınır şartlarının birinde özdeğer parametresi bulunan ve tek noktada süreksiz olan Sturm-Liouville probleminin özdeğerlerini incelemiştir. Biz de bu çalışmada [-1,1] aralığının $h_1$ ve $h_2$ gibi iki iç noktasında süreksiz olan, katsayıları sonlu

$$\ell u := \frac{1}{r(x)}\left\{-\left(p(x)u'\right)' + q(x)u\right\} = \lambda u,\ x \in [-1,h_1) \cup (h_1,h_2) \cup (h_2,1] \qquad (1)$$

diferansiyel denkleminden,

$$u(-1) = 0, \qquad (2)$$

$$(\lambda\alpha_1 + \beta_1)u(1) - (\lambda\alpha_2 + \beta_2)u'(1) = 0 \qquad (3)$$

sınır koşullarından ve x=$h_1$, x=$h_2$ süreksizlik noktalarındaki

$$\gamma_1 u(h_1 - 0) = \delta_1 u(h_1 + 0), \tag{4}$$

$$\gamma_2 u'(h_1 - 0) = \delta_2 u'(h_1 + 0), \tag{5}$$

$$\gamma_3 u(h_2 - 0) = \delta_3 u(h_2 + 0), \tag{6}$$

$$\gamma_4 u'(h_2 - 0) = \delta_4 u'(h_2 + 0), \tag{7}$$

geçiş koşullarından oluşan sınır-değer problemi incelendi. Burada $\lambda$ kompleks parametredir; $\alpha_i, \beta_i$ ($i=1,2$), $\gamma_j, \delta_j$ ($j=1,2,3,4$) reel sayılardır ve $|\beta_1|+|\beta_2|\neq 0$, $|\gamma_j|+|\delta_j|\neq 0$ ($j=1,2,3,4$) doğal koşullarını sağlıyorlar; $r(x)$, $p(x)$, $p'(x)$, $q(x)$ ise [-1,$h_1$), ($h_1,h_2$) ve ($h_2$,1] aralıklarının her birinde sürekli olan x=$h_1$, x=$h_2$ noktalarında sonlu sol ve sağ limit değerleri mevcut olan reel değerli fonksiyonlardır. Ayrıca her $x \in [-1,h_1) \cup (h_1, h_2) \cup (h_2, 1]$ için $r(x) > 0, p(x) > 0$ olduğunu kabul edeceğiz. Bu problemin [16] anlamında kendine eşlenik olması için, $\rho := \alpha_1 \beta_2 - \alpha_2 \beta_1 > 0$ şartının da sağlandığını kabul edeceğiz. Walter'in [16] makalesinde olduğu gibi, eğer (1)-(7) sınır-değer problemi herhangi bir Hilbert uzayında kendine eşlenik bir operatör için özdeğer problemine indirgenebilirse, o halde bu probleme kendine eşlenik problem diyeceğiz. (1)–(7) probleminin bazı özel halleri [1, 6, 7, 10, 11] incelenmiştir.

Matematiksel fiziğin bazı problemlerinde zaman değişkenine göre kısmi türev sadece diferansiyel denklemde değil aynı zamanda sınır koşularında da ortaya çıkmaktadır. Böyle problemlere uygun olan spektral problemlerde özdeğer parametresi sadece diferansiyel denklemde değil aynı zamanda sınır koşullarında da bulunmaktadır [4, 13]. (4)-(7) biçimindeki geçiş koşullarında ise farklı fiziksel ve mekanik özellikleri bulunan cisimler arasındaki ısı ve madde iletimi veya başka geçiş süreçlerinde ortaya çıkmaktadır [5, 6, 9-11, 15].

## SINIR-DEĞER PROBLEMİNİN UYGUN HILBERT UZAYINDA ÖZDEĞER PROBLEMİ OLARAK İFADESİ

Eğer,

$$\begin{cases} (u)(1) := \beta_1 u(1) - \beta_2 u'(1) \\ (u)'(1) := \alpha_1 u(1) - \alpha_2 u'(1) \end{cases} \tag{8}$$

gösterimlerinden yararlanırsak, kolay bir şekilde $u, v \in C^1[-1,1]$ için,

$$\rho[u(1)v'(1) - u'(1)v(1)] = (u)_1 (v)'_1 - (u)'_1 v(1)_1 \tag{9}$$

olduğunu gösterebiliriz. Şimdi iki bileşenli

$$T := \begin{pmatrix} T_1(x) \\ T_2 \end{pmatrix}, \; T_1(x) \in L_2[-1,1], T_2 \in \mathbb{C} \;\; ; \;\; G := \begin{pmatrix} G_1(x) \\ G_2 \end{pmatrix}, \; G_1(x) \in L_2[-1,1], G_2 \in \mathbb{C} \;\; ;$$

elemanlarının $L_2[-1,1] \oplus \mathbb{C}$ lineer uzayında iki $T, G \in L_2[-1,1] \oplus \mathbb{C}$ elemanlarının iç çarpımını,

$$<T,G>_{p,r,\rho} := \int_{-1}^{h_1} T_1(x)\overline{G_1(x)}r(x)dx + \int_{h_1}^{h_2} T_1(x)\overline{G_1(x)}r(x)dx + \int_{h_2}^{1} T_1(x)\overline{G_1(x)}r(x)dx + \frac{p(1)}{\rho}T_2\overline{G_2}$$

formülü ile tanımlayalım. O halde

$H_{p,r,\rho} := (L_2[-1,1] \oplus \mathbb{C}, <\bullet,\bullet>_{p,r,\rho})$ iç çarpım uzayının bir Hilbert uzayı olacağı açıktır. Bu uzayda tanım bölgesi

D(K)={ $T \in H_{p,r,\rho}$ | $T_1, T_1'$ fonksiyonlarının her biri [-1, $h_1$), ($h_1, h_2$) ve ($h_2, 1$] aralıklarının her birinde mutlak süreklidir; $T_1(h_1 \pm 0)$, $T_1'(h_1 \pm 0)$ $T_1(h_2 \pm 0)$, $T_2'(h_2 \pm 0)$ sonlu limit değerleri mevcuttur. $T_1(-1) = 0$, $\gamma_1 T_1(h_1 - 0) = \delta_1 T_1(h_1 + 0)$, $\gamma_2 T_1'(h_1 - 0) = \delta_2 T_1'(h_1 + 0)$, $\gamma_3 T_1(h_2 - 0) = \delta_3 T_1(h_2 + 0)$, $\gamma_4 T_1'(h_2 - 0) = \delta_4 T_1'(h_2 + 0)$ ve $T_2 = (T_1)_1'$ } (10)

olan $K : H_{p,r,\rho} \to H_{p,r,\rho}$ operatörünü,

$$K \begin{pmatrix} T_1(x) \\ (T_1)_1' \end{pmatrix} := \begin{pmatrix} \ell T_1 \\ -(T_1)_1 \end{pmatrix} \tag{11}$$

eşitliği ile tanımlayalım. O halde (1)-(7) sınır-değer problemi,

$$KU = \lambda U \quad \left( U := \begin{pmatrix} u(x) \\ (u)_1' \end{pmatrix} \in D(K) \right) \tag{12}$$

operatör-denklem biçiminde yazılabilir. Böylece (1)-(7) sınır-değer problemini bir Hilbert uzayında tanımlı olan bir lineer operatör için özdeğer problemine indirgemiş olduk.

**2.1 Lemma.** Eğer $\delta_1 \delta_2 p(h_1 - 0) = \gamma_1 \gamma_2 p(h_1 + 0)$ ve $\delta_3 \delta_4 p(h_2 - 0) = \gamma_3 \gamma_4 p(h_2 + 0)$ şartları sağlanıyorsa $K$ operatörü simetriktir.

**İspat.** $T, G \in D(K)$ herhangi iki eleman olmak üzere iyi-bilinen Lagrange formülünü [8] uygularsak,

$$< KT, G >_{p,r,\rho} := \int_{-1}^{1} (\ell T_1)(x) \overline{G_1(x)} r(x) dx + \frac{p(1)}{\rho} ((-T_1)_1)\overline{(G_1)_1'} = \int_{-1}^{h_1} T_1(x) \overline{(\ell G_1)(x)} r(x) dx$$

$$+ \int_{h_1}^{h_2} T_1(x) \overline{(\ell G_1)(x)} r(x) dx + \int_{h_2}^{1} T_1(x) \overline{(\ell G_1)(x)} r(x) dx + p(h_1 - 0) W(T_1, G_1; h_1 - 0) - p(-1) W(T_1, G_1; -1)$$

$$+ p(h_2 - 0) W(T_1, G_1; h_2 - 0) - p(h_1 + 0) W(T_1, G_1; h_1 + 0) + p(1) W(T_1, G_1; 1)$$

$$- p(h_2 + 0) W(T_1, G_1; h_2 + 0) - \frac{p(1)}{\rho} (T_1)_1 \overline{(G_1)_1'}$$

$$= \left[ \langle T, G \rangle_{p,r,\rho} - \frac{p(1)}{\rho} \overline{(G_1)_1} (T_1)_1' \right] + \left[ p(h_1 - 0) W(T_1, G_1; h_1 - 0) - p(h_1 + 0) \right.$$

$$\times W(T_1, G_1; h_1 + 0) \bigg] + \bigg[ p(h_2 - 0)$$

$$\times W(T_1, G_1; h_2 - 0) - p(h_2 + 0) W(T_1, G_1; h_2 + 0) \bigg] - \frac{p(1)}{\rho} [(T_1)_1 \overline{(G_1)_1'} - (T_1)_1' \overline{(G_1)_1}] =$$

$$\left[ \int_{-1}^{h_1} T_1(x) \overline{(\ell G_1)(x)} r(x) dx + \int_{h_1}^{h_2} T_1(x) \overline{(\ell G_1)(x)} r(x) dx + \int_{h_2}^{1} T_1(x) \overline{(\ell G_1)(x)} r(x) dx + \frac{p(1)}{\rho} (T_1)_1' \overline{((-G_1)_1)} \right]$$

$$- p(-1) W(T_1, \overline{G_1}; -1) + p(h_1 - 0) W(T_1, \overline{G_1}; h_1 - 0)$$

$$- p(h_1 + 0) W(T_1, \overline{G_1}; h_1 + 0) + p(h_2 - 0) W(T_1, \overline{G_1}; h_2 - 0) - p(h_2 + 0)$$

$$\times W(T_1, \overline{G_1}; h_2 + 0) + p(1) \left[ W(T_1, \overline{G_1}; 1) - \frac{1}{\rho} \left( (T_1)_1 \overline{(G_1)_1'} - (T_1)_1' \overline{(G_1)_1} \right) \right] \tag{13}$$

eşitliğini buluruz. Burada;

$$W(T_1, G_1; x) := T_1(x)G_1'(x) - T_1(x)'G_1(x) \tag{14}$$

ile $T_1$, $G_1$ fonksiyonlarının Wronkskiyeni gösterilmiştir. $T_1(x)$ ve $\overline{G_1}(x)$ fonksiyonları (2) sınır koşulunu sağladıkları için,

$$W(T_1, \overline{G_1}; -1) = 0 \tag{15}$$

eşitliği sağlanır. $T_1$, $\overline{G_1}$ fonksiyonlarının (4)-(7) geçiş koşularını sağladığını ve lemmanın koşulunu dikkate alırsak,

$$p(h_2 - 0)W(T_1, \overline{G_1}; h_2 - 0) = p(h_2 - 0)\left[T_1(h_2 - 0)\overline{G_1}'(h_2 - 0) - T_1'(h_2 - 0)\overline{G_1}(h_2 - 0)\right]$$

$$= \frac{\gamma_3 \gamma_4}{\delta_3 \delta_4} p(h_2 + 0)\left\{\left(\frac{\delta_3}{\gamma_3}T_1(h_2 + 0)\right)\left(\frac{\delta_4}{\gamma_4}\overline{G_1}'(h_2 + 0)\right) - \left(\frac{\delta_3}{\gamma_3}T_1'(h_2 + 0)\right)\left(\frac{\delta_4}{\gamma_4}\overline{G_1}(h_2 + 0)\right)\right\}$$

$$-\left(\frac{\delta_3}{\gamma_3}T_1'(h_2 + 0)\right)\left(\frac{\delta_4}{\gamma_4}\overline{G_1}(h_2 + 0)\right)\right\} = p(h_2 + 0)W(T_1, \overline{G_1}; h_2 + 0) \tag{16}$$

ve benzer şekilde

$$p(h_1 - 0)W(T_1, \overline{G_1}; h_1 - 0) = p(h_1 + 0)W(T_1, \overline{G_1}; h_1 + 0) \tag{17}$$

elde ederiz. O halde (15), (16) ve (17) eşitsizliklerini (13) de yerine yazarak ve (9) eşitliğini de göz önüne alarak, arzu edilen $<KT, G>_{p,r,\rho} = \langle T, KG \rangle_{p,r,\rho}$ eşitliğini, yani $K$ operatörünün simetrik olduğunu elde ederiz.

**2.1 Sonuç.** (1)-(7) sınır-değer probleminin bütün özdeğerleri reeldir.

Not: $p(x)$, $q(x)$ ve $r(x)$ reel değerli fonksiyonlar, (2)-(7) koşullarının katsayıları reel sayılar ve bütün özdeğerler reel olduğu için (1)-(7) probleminin bütün özfonksiyonlarını reel değerli fonksiyonlar olarak kabul edebiliriz.

**2.2 Sonuç.** $\lambda_1$ ve $\lambda_2$ (1)-(7) probleminin herhangi iki farklı özdeğeri, $u_1(x)$ ve $u_2(x)$ bunlara karşılık gelen özfonksiyonlar ise

$$\int_{-1}^{1} u_1(x)u_2(x)r(x)dx = \frac{-p(1)}{\rho}(u_1)_1'(u_2)_1' \tag{18}$$

eşitliği sağlanır.

### $K$ OPERATÖRÜNÜN REZOLVENTİ

Bu kesimde özdeğer olmayan her $\lambda \in$ sayısının $K$ operatörünün regüler değeri olduğunu göstereceğiz ve ayrıca, $R(\lambda, K) := (K - \lambda I)^{-1}$ rezolvent operatörünü inceleyeceğiz. Burada $I$ birim matristir. $T \in H_{p,r,\rho}$ keyfi elemanı için

$$(K - \lambda I)U = T \tag{19}$$

operatör denklemini onunla eşdeğer, homojen olmayan

$$\frac{1}{r(x)}\{-(p(x)U_1')' + q(x)U_1\} - \lambda U_1 = T_1(x), x \in [-1, h_1) \cup (h_1, h_2) \cup (h_2, 1] \tag{20}$$

$$U_1(-1) = 0 \tag{21}$$

$$\left(\beta_1 U_1(1) - \beta_2 U_1'(1)\right) + \lambda\left(\alpha_1 U_1(1) - \alpha_2 U_1'(1)\right) = T_2 \tag{22}$$

$$\gamma_1 U_1(h_1 - 0) - \delta_1 U_1(h_1 + 0) = 0 \tag{23}$$

$$\gamma_2 U_1'(h_1 - 0) - \delta_2 U_1'(h_1 + 0) = 0 \tag{24}$$

$$\gamma_3 U_1(h_2 - 0) - \delta_3 U_1(h_2 + 0) = 0 \tag{25}$$

$$\gamma_4 U_1'(h_2 - 0) - \delta_4 U_1'(h_2 + 0) = 0 \tag{26}$$

sınır-değer problemi seklinde yazalım. Öncelikle aşağıdaki lemmayı verelim.

**3.1 Lemma.** Herhangi $[d_1, d_2]$ aralığında tanımlı ve reel değerli $r(x) \neq 0, p(x) \neq 0$ ve $q(x)$ fonksiyonları verilsin. Eğer $r(x)$ ve $q(x)$ fonksiyonları bu aralıkta sürekli, $p(x)$ ise sürekli diferansiyellenebilir ise, o halde her tam $f(\lambda)$ ve $g(\lambda)$ fonksiyonları için

$$\frac{1}{r(x)}\{-(p(x)u')' + q(x)u\} = \lambda u, \quad x \in [d_1, d_2] \tag{27}$$

diferansiyel denkleminin

$$u(d_i) = f(\lambda), u'(d_i) = g(\lambda) \ (i = 1 \text{ veya } 2) \tag{28}$$

başlangıç koşullarını sağlayan $u(x, \lambda)$ çözümü mevcuttur ve bu çözüm fonksiyonu her $x \in [d_1, d_2]$ değeri için $\lambda$ değişkeninin tam fonksiyonudur. Bu lemma Titchmarsh'ın [14] kitabındaki Teorem 1.5'in ispatındaki yöntemle benzer şekilde ispat edilir.

Şimdi bu lemmadan faydalanarak (1) diferansiyel denkleminin $\varphi(x, \lambda)$ ve $\chi(x, \lambda)$ gibi iki tane çözümünü tanımlayacağız. $[-1, h_1]$ aralığında (1) diferansiyel denkleminin $u(-1) = 0$, $u'(-1) = 1$ başlangıç koşullarını sağlayan çözümünü $\varphi_1(x, \lambda)$ ile gösterelim. $\varphi_1(x, \lambda)$ fonksiyonu tanımlandıktan sonra $[h_1, h_2]$ aralığında (1) diferansiyel denkleminin

$$u(h_1) = \frac{\gamma_1}{\delta_1}\varphi_1(h_1 - 0, \lambda), \quad u'(h_1) = \frac{\gamma_2}{\delta_2}\varphi_1'(h_1 - 0, \lambda) \tag{29}$$

başlangıç koşullarını sağlayan çözümünü tanımlayabiliriz. Bu çözümü $\varphi_2(x, \lambda)$ ile gösterelim. Benzer şekilde $[h_2, 1]$ aralığında (1) diferansiyel denkleminin

$$u(h_2) = \frac{\gamma_3}{\delta_3}\varphi_2(h_2 - 0, \lambda), \quad u'(h_2) = \frac{\gamma_4}{\delta_4}\varphi_2'(h_2 - 0, \lambda) \tag{30}$$

başlangıç koşullarını sağlayan çözümünü $\varphi_3(x, \lambda)$ ile gösterelim. Yine benzer şekilde, $[h_2, 1]$ aralığında (1) diferansiyel denkleminin

$$u(1) = \alpha_2 \lambda + \beta_2, u'(1) = \alpha_1 \lambda + \beta_1 \tag{31}$$

başlangıç koşullarını sağlayan çözümünü $\chi_3(x, \lambda)$ ile göstererek, bu çözümü tanımladıktan sonra $[h_1, h_2]$ aralığında (1) diferansiyel denkleminin

$$u(h_2) = \frac{\delta_3}{\gamma_3}\chi_3(h_2 - 0, \lambda), u'(h_2) = \frac{\delta_4}{\gamma_4}\chi_3'(h_2 - 0, \lambda) \tag{32}$$

başlangıç koşullarını sağlayan çözümünü

$\chi_2(x,\lambda)$ ile göstererek, bu çözümü tanımladıktan sonra $[-1,h_1]$ aralığında (1) diferansiyel denkleminin

$$u(h_1) = \frac{\delta_1}{\gamma_1}\chi_2(h_1-0,\lambda),\ u'(h_1) = \frac{\delta_2}{\gamma_2}\chi_2'(h_1-0,\lambda) \qquad (33)$$

başlangıç koşullarını sağlayan çözümünü $\chi_1(x,\lambda)$ ile gösterelim. 3.1 Lemma gereği $\varphi_i(x,\lambda), \chi_i(x,\lambda)$ ($i=1, 2, 3$) fonksiyonları $\lambda$ nın tam fonksiyonlarıdır. Bu fonksiyonların tanımları gereği

$$\chi(x,\lambda) = \begin{cases} \chi_1(x,\lambda), x \in [-1,h_1], \\ \chi_2(x,\lambda), x \in [h_1,h_2], \\ \chi_3(x,\lambda), x \in [h_2,1] \end{cases}$$

eşitlikleri ile tanımlı $\varphi$ ve $\chi$ fonksiyonları $[-1,h_1) \cup (h_1,h_2) \cup (h_2,1]$ de (1) denklemini ve (4)-(7) geçiş koşullarını sağlayacaklardır. Ayrıca $\varphi(x,\lambda)$ çözümü (2) sınır koşulunu, $\chi(x,\lambda)$ ise (3) sınır koşulunu sağlayacaktır. Aşağıdaki;

$\omega_i := W_\lambda(\varphi_i,\chi_i;x)$ (i=1,2)

$$\omega(x,\lambda) := W_\lambda(\varphi_i,\chi_i;x) = \begin{cases} \omega_1(x,\lambda), x \in [-1,h_1), \\ \omega_2(x,\lambda), x \in (h_1,h_2), \\ \omega_3(x,\lambda), x \in (h_2,1] \end{cases}$$

gösterimlerinden faydalanacağız.

**3.2 Lemma.** Özdeğer olmayan her $\lambda \in$ ve her $x \in [-1,h_1) \cup (h_1,h_2) \cup (h_2,1]$ için $\omega(x,\lambda) \neq 0$ dir.

**3.1 Sonuç.** Özdeğer olmayan her $\lambda \in$ için $\varphi_1(x,\lambda), \chi_1(x,\lambda)$ fonksiyonları $[-1,h_1]$ aralığında, $\varphi_2(x,\lambda), \chi_2(x,\lambda)$ fonksiyonları $[h_1,h_2]$ aralığında, $\varphi_3(x,\lambda), \chi_3(x,\lambda)$ fonksiyonları ise $[h_2,1]$ aralığında lineer bağımsızdırlar.

3.1 Sonuç gereği özdeğer olmayan her $\lambda \in$ için (1) diferansiyel denkleminin genel çözümünü

$$u(x,\lambda) = \begin{cases} C_1\varphi_1(x,\lambda) + D_1\chi_1(x,\lambda), x \in [-1,h_1), \\ C_2\varphi_2(x,\lambda) + D_2\chi_2(x,\lambda), x \in (h_1,h_2), \\ C_3\varphi_3(x,\lambda) + D_3\chi_3(x,\lambda), x \in (h_2,1] \end{cases}$$

biçiminde ifade edebiliriz. Burada $C_i$, $D_i$ ($i=1,2,3$) keyfi sabitlerdir. O halde sabitin varyasyonu yöntemini [8] uygulayarak (20) homojen olmayan denkleminin genel çözümünü $x \in [-1,h_1)$ için

$$U_1(x,\lambda) = \chi_1(x,\lambda)\int_{-1}^{x}\frac{\varphi_1(y,\lambda)}{\omega_1(y,\lambda)}T_1(y)dy + \varphi_1(x,\lambda)\int_{x}^{h_1}\frac{\chi_1(y,\lambda)}{\omega_1(y,\lambda)}T_1(y)dy + C_1\varphi_1(x,\lambda) + D_1\chi_1(x,\lambda) \qquad (34)$$

biçiminde, $x \in (h_1,h_2)$ için

$$U_1(x,\lambda) = \chi_2(x,\lambda)\int_{h_1}^{x}\frac{\varphi_2(y,\lambda)}{\omega_2(y,\lambda)}T_1(y)dy + \varphi_2(x,\lambda)\int_{x}^{h_2}\frac{\chi_2(y,\lambda)}{\omega_2(y,\lambda)}T_1(y)dy + C_2\varphi_2(x,\lambda) + D_2\chi_2(x,\lambda) \quad (35)$$

biçiminde, $x \in (h_2, 1]$ için ise

$$U_1(x,\lambda) = \chi_3(x,\lambda)\int_{h_2}^{x}\frac{\varphi_3(y,\lambda)}{\omega_3(y,\lambda)}T_1(y)dy + \varphi_3(x,\lambda)\int_{x}^{1}\frac{\chi_3(y,\lambda)}{\omega_3(y,\lambda)}T_1(y)dy + C_3\varphi_3(x,\lambda) + D_3\chi_3(x,\lambda) \quad (36)$$

biçiminde ifade edebiliriz. (20) diferansiyel denkleminin (34)-(36) eşitlikleri ile verilmiş genel çözümünü (21)-(26) koşularında yerine yazarsak $C_i$, $D_i$ sabitlerini bulabiliriz. (34) ifadesini (21) sınır koşulunda yerine yazarsak $D_1\chi(-1,\lambda) = 0$ eşitliğini elde ederiz. $\lambda$ özdeğer olmadığı için $\chi(-1,\lambda) \neq 0$ dır. Dolayısıyla $D_1 = 0$ olur. (36) ifadesini (22) sınır koşulunda yerine yazarsak, $C_3 = \dfrac{T_2}{\omega_3(1,\lambda)}$ buluruz. $D_1$ ve $C_3$ için bulduğumuz değerleri de göz önüne alarak (34)-(36) ifadelerini (23)-(26) geçiş koşularında yazarsak; $C_1, C_2, D_2, D_3$ değerlerini bulmak için aşağıdaki lineer denklem sistemini elde ederiz:

$$\gamma_1\varphi_1(h_1,\lambda)C_1 - \delta_1\chi_2(h_1,\lambda)D_2 = \gamma_1\chi_1(h_1,\lambda)$$
$$\times\int_{-1}^{h_1}\frac{\varphi_1(y,\lambda)}{\omega_1(y,\lambda)}T_1(y)dy + \delta_1\varphi_2(h_1,\lambda)\int_{h_1}^{h_2}\frac{\chi_2(y,\lambda)}{\omega_2(y,\lambda)}T_1(y)dy + \delta_1 C_2\varphi_2(h_1,\lambda),$$

$$\gamma_2\varphi_1'(h_1,\lambda)C_1 - \delta_2\chi_2'(h_1,\lambda)D_2 = -\gamma_2\chi_1'(h_1,\lambda)$$
$$\times\int_{-1}^{h_1}\frac{\varphi_1(y,\lambda)}{\omega_1(y,\lambda)}T_1(y)dy + \delta_2\varphi_2'(h_1,\lambda)\int_{h_1}^{h_2}\frac{\chi_2(y,\lambda)}{\omega_2(y,\lambda)}T_1(y)dy + \delta_2 C_2\varphi_2'(h_1,\lambda),$$

$$\gamma_3\varphi_2(h_2,\lambda)C_2 - \delta_3\chi_3(h_2,\lambda)D_3 = -\gamma_3\chi_2(h_2,\lambda)$$
$$\times\int_{h_1}^{h_2}\frac{\varphi_2(y,\lambda)}{\omega_2(y,\lambda)}T_1(y)dy + \delta_3\varphi_3(h_2,\lambda)\int_{h_2}^{1}\frac{\chi_3(y,\lambda)}{\omega_3(y,\lambda)}T_1(y)dy + \frac{\delta_3 T_2\varphi_3(h_2,\lambda)}{\omega_3(1,\lambda)} + \gamma_3 D_2\chi_2(h_2,\lambda),$$

$$\gamma_4\varphi_2'(h_2,\lambda)C_2 - \delta_4\chi_3'(h_2,\lambda)D_3 = -\gamma_4\chi_2'(h_2,\lambda)$$
$$\times\int_{h_1}^{h_2}\frac{\varphi_2(y,\lambda)}{\omega_2(y,\lambda)}T_1(y)dy + \delta_4\varphi_3'(h_2,\lambda)\int_{h_2}^{1}\frac{\chi_3(y,\lambda)}{\omega_3(y,\lambda)}T_1(y)dy + \frac{\delta_4 T_2\varphi_3'(h_2,\lambda)}{\omega_3(1,\lambda)} + \gamma_4 D_2\chi_2'(h_2,\lambda).$$

Bu sistemin determinantı $-\delta_1\delta_2\delta_3\delta_4\omega_2(h_1,\lambda)\omega_3(h_2,\lambda) \neq 0$ olduğu için bir tek çözümü bulunur. $\varphi_i(x,\lambda), \chi_i(x,\lambda)$ fonksiyonlarının tanımlarından yararlanarak sonuncu denklem sisteminden

$$C_1 = \int_{h_1}^{h_2}\frac{\chi_2(y,\lambda)}{\omega_2(y,\lambda)}T_1(y)dy + \int_{h_2}^{1}\frac{\chi_3(y,\lambda)}{\omega_3(y,\lambda)}T_1(y)dy + \frac{T_2}{\omega_3(y,\lambda)},$$

$$C_2 = \int_{h_2}^{1}\frac{\chi_3(y,\lambda)}{\omega_3(y,\lambda)}T_1(y)dy + \frac{T_2}{\omega_3(y,\lambda)}, \quad D_2 = \int_{-1}^{h_1}\frac{\varphi_1(y,\lambda)}{\omega_1(y,\lambda)}T_1(y)dy, \quad D_3 = \int_{h_1}^{h_2}\frac{\varphi_2(y,\lambda)}{\omega_2(y,\lambda)}T_1(y)dy$$

elde edilir. $C_i$, $D_i$ sabitleri için bulduğumuz değerleri (34)-(36) ifadelerinde yerine yazarsak ve gerekli düzenlemeleri yaparsak, (20)-(26) probleminin çözümü için tüm $[-1, h_1) \cup (h_1, h_2) \cup (h_2, 1]$ aralığında

$$U_1 = \chi(x,\lambda)\int_{-1}^{x}\frac{\varphi(y,\lambda)}{\omega(y,\lambda)}T_1(y)dy + \varphi(x,\lambda)\int_{-1}^{x}\frac{\chi(y,\lambda)}{\omega(y,\lambda)}T_1(y)dy + \frac{T_2\varphi(x,\lambda)}{\omega_3(1,\lambda)}$$

formülünü elde ederiz.

3.**1. Teorem**. Özdeğer olmayan her $\lambda \in$ sayısı (10), (11) eşitlikleri ile tanımlı olan $K$ operatörünün regüler değeridir. Ayrıca $R(\lambda,K): H_{p,r,\rho} \to H_{p,r,\rho}$ rezolvent operatörü kompakt operatördür.

**İspat.**

$$G_1(x,y;\lambda) = \begin{cases} \dfrac{\chi(x,\lambda)\varphi(y,\lambda)}{\omega(y,\lambda)}, & -1 \le y \le x \le 1; \\ & x,y \ne h_i (i=1,2), \\ \dfrac{\varphi(x,\lambda)\chi(y,\lambda)}{\omega(y,\lambda)} & -1 \le x \le y \le 1; \\ & x,y \ne h_i (i=1,2) \end{cases}$$ gösteriminden yararlanarak sonuncu formülü

$U_1(x,\lambda) = \int_{-1}^{1} G_1(x,y;\lambda)T_1(y)dy + \dfrac{T_2}{\omega(1,\lambda)}\varphi(x,\lambda)$ biçiminde ifade edebiliriz. Buradan $R(\lambda,K)$ rezolvent operatörü için

$$R(\lambda,K)T = \begin{pmatrix} \int_{-1}^{1} G_1(x,y;\lambda)T_1(y)dy + \dfrac{T_2}{\omega(1,\lambda)}\varphi(x,\lambda) \\ \int_{-1}^{1} \left(G_1(\bullet,y;\lambda)\right)_1' T_1(y)dy + \dfrac{T_2}{\omega(1,\lambda)}\left(\varphi(\bullet,\lambda)\right)_1' \end{pmatrix}$$ formülü elde edilir.

$B_\lambda: L_2[-1,1] \to L_2[-1,1], B_\lambda: H_{p,r,\rho} \to H_{p,r,\rho}$ ve $S_\lambda: H_{p,r,\rho} \to H_{p,r,\rho}$ operatörlerini

$$B_\lambda T_1 := \int_{-1}^{1} G_1(x,y;\lambda)T_1(y)dy, \quad B_\lambda T := \begin{pmatrix} B_\lambda T_1 \\ (B_\lambda T_1)_1' \end{pmatrix}, \quad S_\lambda T := \begin{pmatrix} \dfrac{T_2}{\omega(1,\lambda)}\varphi(x,\lambda) \\ \dfrac{T_2}{\omega(1,\lambda)}\left(\varphi(\bullet,\lambda)\right)_1' \end{pmatrix}$$

eşitlikleri ile tanımlarsak, $R(\lambda,K)$ rezolvent operatörünü $R(\lambda,K) = B_\lambda + S_\lambda$ biçiminde ifade edebiliriz. $B_\lambda$ operatörü $L_2[-1,1]$ Hilbert uzayında kompakt olduğu için, $B_\lambda$ operatörü $H_{p,r,\rho}$ Hilbert uzayında kompakttır [8]. $S_\lambda$ operatörünün $H_{p,r,\rho}$ Hilbert uzayında kompakt olduğu açıktır. Dolayısıyla özdeğer olmayan her $\lambda \in$ için $R(\lambda,K)$ operatörü de $H_{p,r,\rho}$ uzayında kompakttır.

## ÖZFONKSİYONLAR SİSTEMİNİN SERİSİNE AÇILIM

**4.1 Teorem.** (10), (11) eşitlikleri ile tanımlı $K$ operatörü $H_{p,r,\rho}$ Hilbert uzayında kendine eşleniktir.

Sonuç olarak, 3.1 Teorem, 4.1 Teorem ve iyi bilinen Hilbert-Schmidt Teoremi [12] gereği aşağıdaki teorem elde edilir [2].

**4.2 Teorem**. $H_{p,r,\rho}$ Hilbert uzayında (10), (11) eşitlikleri ile tanımlı $K$ operatörünün sayılabilir sayıda reel özdeğeri mevcuttur, her özdeğerin cebirsel katı sonludur, özdeğerler dizisi alttan sınırlıdır ve sonlu yığılma noktası yoktur. Her özdeğer cebirsel katı sayıda yazılmak kaydı ile, özdeğerler dizisini $\lambda_1 \leq \lambda_2 \leq ...$ biçiminde sıralayarak, uygun normlandırılmış özelementler

$$\varphi_n := \begin{pmatrix} \phi_n(x) \\ (\phi_n)'_1 \end{pmatrix}, \left( \|\varphi_n\|_{H_{p,r,\rho}} = 1, \, n = 1, 2, ... \right)$$

biçiminde gösterilmek üzere, her $T \in H_{p,r,\rho}$ elemanı için $\sum_{n=1}^{\infty} c_n \varphi_n, \, c_n = \langle T, \varphi_n \rangle_{H_{p,r,\rho}}$ Fourier serisi $H_{p,r,\rho}$ Hilbert uzayında $T$ elemanına yakınsak olacaktır;

$$T = \sum_{n=1}^{\infty} \langle T, \varphi_n \rangle_{H_{p,r,\rho}} \varphi_n . \tag{37}$$

**4.1 Sonuç.** Her $f \in L_2[-1,1]$ fonksiyonu $L_2([-1,1], r)$ Hilbert uzayında (1)-(7) sınır-değer probleminin $\{\phi_n\}, \, n = 1, 2, ...$ özfonksiyonlar sisteminin $f(x) = \sum_{n=1}^{\infty} \left( \int_{-1}^{1} f(y) \phi_n(y) r(y) dy \right) \phi_n(x)$ serisine açılır.

**İspat.** İspat için (37) formülünde $T \in H_{p,r,\rho}$ elemanını özel olarak $T = \begin{pmatrix} f(x) \\ 0 \end{pmatrix}$ almak yeterlidir. (Burada $L_2([-1,1], r)$ ve $L_2[-1,1]$ Hilbert uzaylarının lineer uzaylar olarak aynı olduğuna dikkat etmek gerekir.)

**4.2 Sonuç.** Her $f \in L_2[-1,1]$ fonksiyonu için,

$$\sum_{n=1}^{\infty} [(\phi_n)'_1]^2 = \frac{\rho}{p(1)}, \tag{38}$$

$$\sum_{n=1}^{\infty} (\phi_n)'_1 \phi_n(x) = 0 \tag{39}$$

eşitlikleri sağlanır.

**İspat.** (37) formülünü

$$\begin{pmatrix} T_1(x) \\ T_2 \end{pmatrix} = \begin{pmatrix} \sum_{n=0}^{\infty} \langle T, \varphi_n \rangle_{H_{p,r,\rho}} \phi_n(x) \\ \sum_{n=0}^{\infty} \langle T, \varphi_n \rangle_{H_{p,r,\rho}} (\phi_n)'_1 \end{pmatrix} \tag{40}$$

biçiminde yazalım. Bu formülde özel olarak

$T = \begin{pmatrix} 0 \\ 1 \end{pmatrix}$ alırsak, $\begin{pmatrix} 0 \\ 1 \end{pmatrix} = \begin{pmatrix} \sum_{n=0}^{\infty} \frac{p(1)}{\rho} (\phi_n)'_1 \phi_n(x) \\ \sum_{n=0}^{\infty} \frac{p(1)}{\rho} [(\phi_n)'_1]^2 \end{pmatrix}$ eşitliğini, yani (38) ve (39) eşitliklerini elde ederiz.

**4.3 Sonuç**. Her $f \in L_2[-1,1]$ için $\sum_{n=0}^{\infty}\left(\int_{-1}^{1} f(y)\phi_n(y)dy\right)(\phi_n)'_1 = 0$ eşitliği sağlanır.

**İspat**. İspat için (40) formülünü $T = \begin{pmatrix} f(x) \\ 0 \end{pmatrix}$ elemanı için yazmak yeterlidir.

## Kaynaklar